\documentclass[12pt]{article}

\usepackage{amsmath, epsfig, cite}
\usepackage{amssymb}
\usepackage{amsfonts}
\usepackage{latexsym}
\usepackage{amsthm}

\makeatletter
\renewcommand{\@seccntformat}[1]{{\csname the#1\endcsname}.\hspace{.5em}}
\makeatother

\makeatletter
\def\@xnthm#1#2{%
  \let\@tempa\relax
  \@xp\@ifdefinable\csname #2 .\endcsname{%
    \global\@xp\let\csname end#2\endcsname\@endtheorem
    \ifx *#1
      \edef\@tempa##1{%
        \gdef\@xp\@nx\csname#2\endcsname{%
          \@nx\@thm{\@xp\@nx\csname th@\the\thm@style\endcsname}%
            {}{##1}}}%
    \else 
      \def\@tempa{\@oparg{\@ynthm{#2}}[]}%
    \fi
  }%
  \@tempa
}
\def\thmhead@plain#1#2#3{%
  \thmname{#1}\thmnumber{\@ifnotempty{#1}{ }#2}%
  \thmnote{ {\the\thm@notefont#3}}}
\let\thmhead\thmhead@plain
\def\swappedhead#1#2#3{%
  \thmnumber{#2}\thmname{\@ifnotempty{#2}{. }#1}%
  \thmnote{ {\the\thm@notefont#3}}}
\makeatother

\newtheorem{thm}{Theorem}

\newtheorem{lem}[thm]{Lemma}

\newcommand{\pf}{\noindent{\it Proof.} }

\def\k{{\bf k}}
\def\z{{\bf z}}
\def\a{{\bf a}}
\def\b{{\bf b}}
\def\e{{\bf e}}
\def\0{{\bf 0}}
\def\n{{\bf n}}

\begin{document}

\begin{center}
{\Large\bf Bijective proofs of Gould-Mohanty's and\\[3pt]Raney-Mohanty's
identities}
\end{center}

\vskip 2mm \centerline{\large Victor J. W. Guo}


\vskip 0.7cm \noindent{\small{\bf Abstract.} Using the model of
words, we give bijective proofs of Gould-Mohanty's and
Raney-Mohanty's identities, which are respectively multivariable
generalizations of Gould's identity
\begin{align*}
\sum_{k=0}^{n}{x-kz\choose k}{y+kz\choose n-k}=
\sum_{k=0}^{n}{x+\epsilon-kz\choose k}{y-\epsilon+kz\choose n-k}
\end{align*}
and Rothe's identity
\begin{align*}
\sum_{k=0}^{n}\frac{x}{x-kz}{x-kz\choose k}{y+kz\choose n-k}=
{x+y\choose n}.
\end{align*}}



\section{Introduction}
A famous generalization of the binomial
theorem is Abel's identity \cite{Abel}:
\begin{align}
\sum_{k=0}^{n}{n\choose k}x(x-kz)^{k-1}(y+kz)^{n-k}=(x+y)^n, \label{eq:abel-1}
\end{align}
which also has a company identity as follows:
\begin{align}
\sum_{k=0}^{n}{n\choose k}xy(x-kz)^{k-1}(y+kz)^{n-k-1}
=(x+y+nz)(x+y)^{n-1}.
\label{eq:abel-2}
\end{align}
It is not difficult to see that \eqref{eq:abel-1} and \eqref{eq:abel-2} are
respectively limiting cases of the following convolution formulas due to Rothe \cite{Rothe}:
\begin{align}
& \sum_{k=0}^{n}\frac{x}{x-kz}{x-kz\choose k}{y+kz\choose n-k}=
{x+y\choose n}, \label{eq:gou-2} \\[5pt]
& \sum_{k=0}^{n}\frac{xy}{(x-kz)(y-(n-k)z)}{x-kz\choose k}{y-(n-k)z\choose n-k} \nonumber\\
&\phantom{\sum_{k=0}^{n}\frac{x}{x-kz}{x-kz\choose k}{y+kz\choose n-k}}=
\frac{x+y}{x+y-nz}{x+y-nz\choose n}. \label{eq:gou-1}
\end{align}

Gould \cite{Gould,Gould57} reproved \eqref{eq:gou-2} and \eqref{eq:gou-1} and
also obtained the following identity
\begin{align}
\sum_{k=0}^{n}{x-kz\choose k}{y+kz\choose n-k}=
\sum_{k=0}^{n}{x+\epsilon-kz\choose k}{y-\epsilon+kz\choose n-k}. \label{eq:gou-3}
\end{align}
Another proof of \eqref{eq:gou-2} and \eqref{eq:gou-1} was given by
Sprugnoli \cite{Sprugnoli}. It is not difficult to see that
\eqref{eq:gou-1} can be deduced from \eqref{eq:gou-2}. Blackwell and
Dubins \cite{Blackwell} gave a combinatorial proof of Rothe's
identity \eqref{eq:gou-1}, which can also be proved in the model of
lattice paths (using \cite[p.~9]{Mohanty} or
\cite[(1.1)]{Krattenthaler}). Recently, the author \cite{Guo} gives
simple bijective proofs of Gould's identity \eqref{eq:gou-3} and
Rothe's identity \eqref{eq:gou-2} in the model of binary words.

Hurwitz \cite{Hurwitz} established a multivariable generalization of
Abel's identities \eqref{eq:abel-1} and \eqref{eq:abel-2}
(see also \cite{Strehl}). For a curious $q$-analogue of Rothe's identity
\eqref{eq:gou-2}, we refer the reader to \cite{Schlosser} and
references therein.

In order to state a multivariable generalization of Rothe's identities
in the literature, we need first to introduce some notation. Let $m$ be a fixed natural number
throughout the paper.
For $\a=(a_1,\ldots,a_m)\in\mathbb{N}^m$ and
$\b=(b_1,\ldots,b_m)\in\mathbb{C}^m$, set $|\a|=a_1+\cdots+a_m$,
$\a!=a_1!\cdots a_m!$,
$\a+\b=(a_1+b_1,\ldots,a_m+b_m)$, $\a\cdot\b=a_1 b_1+\cdots+a_m b_m$,
and $\b^\a=b_1^{a_1}\cdots b_m^{a_m}$.
For any complex parameter $x$ and $\n=(n_1,\ldots,n_m)\in\mathbb{Z}^m$,
we define the {\it multinomial coefficient} ${x\choose \n}$ by
\begin{align*}
{x\choose \n}=
\begin{cases}
x(x-1)\cdots(x-|\n|+1)/\n!,
&\text{if $\n=(n_1,\ldots,n_m)\in\mathbb{N}^m$,}\\
0, &\text{otherwise.}
\end{cases}
\end{align*}

Using generating functions, Mohanty \cite{Mohanty66} proved the following
multivariable generalization of Rothe's identities \eqref{eq:gou-2} and \eqref{eq:gou-1}:
\begin{align}
& \sum_{\k=\0}^{\n}\frac{x}{x-\k\cdot\z}
{x-\k\cdot\z\choose \k}{y+\k\cdot\z\choose \n-\k}=
{x+y\choose \n}, \label{eq:multi-2} \\[5pt]
& \sum_{\k=\0}^{\n}\frac{xy}{(x-\k\cdot\z)(y-(\n-\k)\cdot\z)}
{x-\k\cdot \z\choose \k}{y-(\n-\k)\cdot\z\choose \n-\k} \nonumber\\
&\hskip 4.5cm{}=\frac{x+y}{x+y-\n\cdot\z}{x+y-\n\cdot\z\choose \n}. \label{eq:multi-1}
\end{align}
However, an important special case of \eqref{eq:multi-1} (where $z_i=i$) was already
contained in the earlier work of Raney \cite{Raney} on a
combinatorial approach to the Lagrange inversion. Hence we would call both
\eqref{eq:multi-2} and \eqref{eq:multi-1} {\it Raney-Mohanty's identities}. Unaware of Mohanty's work,
in 1988 Louck \cite{Louck} proposed a ``conjecture" equivalent to \eqref{eq:multi-1},
which caught the interests of three different people independently and was solved
by them by three different methods:
Paule \cite{Paule} proved \eqref{eq:multi-1} by the Lagrange inversion approach,
Strehl \cite{Strehl} gave a completely combinatorial approach, while Zeng \cite{Zeng}
used mathematical induction.

Moreover, Mohanty and Handa \cite{MH} established the following identity
\begin{align}
\sum_{\k=\0}^{\n}{x+\k\cdot\z\choose \k}{y-\k\cdot\z\choose \n-\k}
=\sum_{\k=\0}^{\n}{x+y-|\k|\choose \n-\k}{|\k|\choose \k}\z^{\k},
\label{eq:mh}
\end{align}
which is a multivariable generalization of Jensen's identity \cite{Gould60}:
\begin{align*}
\sum_{k=0}^{n}{x+kz\choose k}{y-kz\choose n-k}=
\sum_{k=0}^{n}{x+y-k\choose n-k}z^k. 
\end{align*}
It follows immediately from Mohanty-Handa's identity \eqref{eq:mh} that
\begin{align}
\sum_{\k=\0}^{\n}{x-\k\cdot\z\choose \k}{y+\k\cdot\z\choose \n-\k}
=\sum_{\k=\0}^{\n}{x+\epsilon-\k\cdot\z\choose \k}
{y-\epsilon+\k\cdot\z\choose \n-\k}. \label{eq:gmh}
\end{align}
Since \eqref{eq:gmh} is obviously a multivariable generalization of Gould's
identity \eqref{eq:gou-3} and it also follows from one of the generating functions
established by Mohanty in \cite{Mohanty66}, we call \eqref{eq:gmh}
{\it Gould-Mohanty's identity}.

To the knowledge of the author, there are no combinatorial proofs of
Mohanty-Handa's identity \eqref{eq:mh} and Gould-Mohanty's identity
\eqref{eq:gmh}. In this paper, continuing the work of \cite{Guo},
we shall give bijective proofs of Gould-Mohanty's identity and Raney-Mohanty's
identity \eqref{eq:multi-2} in the model of words.


\section{Proof of Gould-Mohanty's identity}\label{sec:gm}
It suffices to prove Gould-Mohanty's identity \eqref{eq:gmh} for the special case:
\begin{align}
\sum_{\k=\0}^{\n}{p-\k\cdot\z\choose \k}{q+\k\cdot\z\choose \n-\k}
=\sum_{\k=\0}^{\n}{p+1-\k\cdot\z\choose \k}
{q-1+\k\cdot\z\choose \n-\k}, \label{eq:gmh-pq}
\end{align}
where $ p,q\in \mathbb{N}$ and $\n,\z\in\mathbb{N}^m$.
Furthermore, we need only to prove that \eqref{eq:gmh-pq} holds for all integers $p\geq\n\cdot\z$
and $q\geq 1$. In this case, each multinomial coefficient in \eqref{eq:gmh-pq} is nonnegative and
therefore has a combinatorial interpretation.

Let $\Gamma=\{a,b_1,\ldots,b_m\}$ denote an alphabet with a grading $||a||=1$
and $||b_i||=z_i+1$ ($1\leq i\leq m$). For a word $w=w_1\cdots w_n\in \Gamma^*$, its {\it length}
$n$ is denoted by $|w|$ and its {\it weight} by
$||w||=||w_1||+\cdots+||w_n||$, and we call the word $w_nw_{n-1}\cdots w_1$
the {\it reverse} of $w$. Let $|w|_{b_i}$ be the number of $b_i$'s
appearing in $w$, and let
$$
\Gamma_{p,\k}:=\{w\in \Gamma^*\colon ||w||=p\
\text{and}\ |w|_{b_i}=k_i,\ i=1,\ldots,m\},
$$
where $\k=(k_1,\ldots,k_m)$. It is easy to see that
$\Gamma_{p,\k}\subseteq \Gamma^{p-\k\cdot\z}$ and
\begin{align}
\#\Gamma_{p,\k}={p-\k\cdot\z\choose \k},
\label{eq:pk}
\end{align}
where $\z=(z_1,\ldots,z_m)$.

Furthermore, let
$$
\Gamma_{p,\k}^{(r)}:=\{w\in \Gamma_{p,\k}\colon
\text{$w$ has a prefix of weight $r$}\}.
$$
For $p,q\geq \n\cdot\z$, an obvious bijection
$$
\Gamma_{p+q,\n}^{(p)}
\longleftrightarrow\biguplus_{\k} \Gamma_{p,\k}\times \Gamma_{q,\n-\k}
$$
leads to
\begin{align}
\#\Gamma_{p+q,\n}^{(p)}
=\sum_{\k}{p-\k\cdot\z\choose \k}{q-(\n-\k)\cdot\z\choose \n-\k}.
\label{eq:pq-p}
\end{align}
Thus, the identity \eqref{eq:gmh-pq} is equivalent to
\begin{equation}\label{eq:pq-num}
\#\Gamma_{p+q+\n\cdot\z,\n}^{(p)}=\#\Gamma_{p+q+\n\cdot\z,\n}^{(p+1)}.
\end{equation}

We need the following simple fact.
\begin{lem}\label{lem:ref}
Let $u,v\in \Gamma^*$ with $||u||,||v||\geq \n\cdot\z+1$, where
$n_i=|u\cdot v|_{b_i}$ {\rm(}$1\leq i\leq m${\rm)}. Then there exist nonempty prefixes $x$
of $u$ and $y$ of $v$ such that $||x||=||y||$.
\end{lem}

\pf Since the proof is easy and very similar to the proof of
\cite[Lemma~1]{Guo}, we omit it here. \qed

Now we can prove \eqref{eq:pq-num} by the following theorem.
\begin{thm}\label{thm:pq}
For all $p\geq \n\cdot\z$ and $q\geq 1$, there is a bijection between
$\Gamma_{p+q+\n\cdot\z,\n}^{(p)}$ and $\Gamma_{p+q+\n\cdot\z,\n}^{(p+1)}$.
\end{thm}
\pf
Suppose that $w=u\cdot v\in \Gamma_{p+q+\n\cdot\z,\n}^{(p)}$,
where $||u||=p$ and $||v||=q+\n\cdot\z$. Applying Lemma~\ref{lem:ref} to $v$ and
the reverse of $u\cdot a$, one sees that $u$ has a suffix $x$ (perhaps empty), i.e.,
$u=u'\cdot x$, and $v$ has a prefix $y$, i.e., $v=y\cdot v'$, such that
$||x||=||y||-1$.
Choosing such $x$ and $y$ with minimal length,
then $w'=u'\cdot \overline y\cdot \overline x\cdot v'\in \Gamma_{p+q+\n\cdot\z,\n}^{(p+1)}$
and $w\mapsto w'$ is a bijection. Here $\overline x$ and $\overline y$
are respectively the reverses of $x$ and $y$.
\qed

In the same manner, we may also give a direct bijection from $\Gamma_{p+q+\n\cdot\z,\n}^{(p)}$
to $\Gamma_{p+q+\n\cdot\z,\n}^{(p+r)}$ for all $p\geq \n\cdot\z$ and $q\geq r\geq 1$.

\section{Proof of Raney-Mohanty's identity}
We again assume that $p\geq \n\cdot\z$ and $q\geq 1$. Moreover, let $z_i\geq 1$ for all $i$.
For each $w\in \Gamma_{p+q+\n\cdot\z,\n}$, let $w=u\cdot v$ denote the
unique factorization with $||u||\geq p$ but as small as possible.
Then we have the following possibilities:
\begin{itemize}
\item If $||u||=p$, then $w\in \Gamma_{p+q+\n\cdot\z,\n}^{(p)}$ and all
these words have been counted in Section~\ref{sec:gm}.

\item If $||u||=p+j$ for some $1\leq j\leq \max\{z_1,\ldots,z_m\}$, then the last letter
of $u$ must a $b_i$ for some $1\leq i\leq m$. Namely, $u=u'\cdot
b_i$ for some $u'\in \Gamma_{p+j-z_i-1,\k-\e_i}$, where
$\e_i=(0,\ldots,1,\ldots,0)\in\mathbb{N}^m$ with the $1$ being in
the $i$-th position. The corresponding $v$ belongs to
$\Gamma_{q+\n\cdot\z-j,\n-\k}$. It is clear that the mapping
$w\mapsto (u',v)$ may be inverted.
\end{itemize}
Hence there is a bijection
\begin{align*}
\Gamma_{p+q+\n\cdot\z,\n}\longleftrightarrow
\Gamma_{p+q+\n\cdot\z,\n}^{(p)}
\biguplus_{i=1}^m\biguplus_{j=1}^{z_i}\biguplus_{\k=\0}^{\n}
\Gamma_{p+j-z_i-1,\k-\e_i}\times \Gamma_{q+\n\cdot\z-j,\n-\k},
\end{align*}
which, together with \eqref{eq:pk} and \eqref{eq:pq-p}, gives the identity
\begin{align}
&\sum_{\k=\0}^{\n}\Bigg({p-\k\cdot\z\choose \k}{q+\k\cdot\z \choose \n-\k}
\nonumber\\[5pt]
&\phantom{\sum_{\k=\0}^{\n}\Bigg(}+\sum_{i=1}^{m}\sum_{j=1}^{z_i}
{p-\k\cdot\z+j-1\choose \k-\e_i}{q+\k\cdot\z-j\choose \n-\k}\Bigg)
={p+q\choose \n}. \label{eq:kmx}
\end{align}
However, by \eqref{eq:gmh}, for all $1\leq i\leq m$ and $1\leq j\leq
z_i$, we have
\begin{equation}\label{eq:kmpink}
\sum_{\k=\0}^{\n}{p-\k\cdot\z+j-1\choose \k-\e_i}{q+\k\cdot\z-j\choose \n-\k}
=\sum_{\k=\0}^{\n}{p-\k\cdot\z-1\choose \k-\e_i}{q+\k\cdot\z\choose \n-\k}.
\end{equation}
Substituting \eqref{eq:kmpink} into \eqref{eq:kmx}, we obtain
\begin{equation}\label{eq:pkroth}
\sum_{\k=\0}^{\n}\left({p-\k\cdot\z\choose \k}
+\sum_{i=1}^m z_{i}{p-\k\cdot\z-1\choose \k-\e_i}\right)
{q+\k\cdot\z \choose \n-\k}={p+q\choose \n}.
\end{equation}
Noticing that
$$
{p-\k\cdot\z-1\choose \k-\e_i}=\frac{k_i}{p-\k\cdot\z}{p-\k\cdot\z\choose \k},
$$
the identity \eqref{eq:pkroth} may be simplified as
\begin{equation*}
\sum_{\k=\0}^{\n}\frac{p}{p-\k\cdot\z}{p-\k\cdot\z\choose \k}
{q+\k\cdot\z\choose \n-\k}={p+q\choose \n},
\end{equation*}
which is Raney-Mohanty's identity \eqref{eq:multi-2}.

For the $m=1$ case, the above bijection also leads
to a double sum extension of the $q$-Chu-Vandermonde
formula (see \cite{Guo}). It is also possible to give
a similar $q$-analogue of \eqref{eq:kmx}. However we omit it here
and leave it to the interested reader.

\section{Some remarks}
We point out that \eqref{eq:multi-1} is a consequence of \eqref{eq:multi-2}, since
the left-hand side of the former may be written as
\begin{align*}
&\hskip -3mm\frac{1}{x+y-\n\cdot\z}\left( \sum_{\k=\0}^{\n}\frac{xy}{x-\k\cdot\z}
{x-\k\cdot \z\choose \k}{y-(\n-\k)\cdot\z\choose \n-\k}\right. \\
&{}\left.+\sum_{\k=\0}^{\n}\frac{xy}{y-(\n-\k)\cdot\z}
{x-\k\cdot \z\choose \k}{y-(\n-\k)\cdot\z\choose \n-\k}\right).
\end{align*}

It is also worth mentioning that Mohanty-Handa's identity \eqref{eq:mh} can
be deduced from Raney-Mohanty's identity \eqref{eq:multi-2}. Indeed, note that
\begin{align*}
\sum_{\k=\0}^{\n}{x+\k\cdot\z\choose \k}{y-\k\cdot\z\choose \n-\k}
&=\sum_{\k=\0}^{\n}\frac{x}{x+\k\cdot\z}
{x+\k\cdot\z\choose \k}{y-\k\cdot\z\choose \n-\k} \nonumber\\
&\quad{}+\sum_{\k=\0}^{\n}\frac{\k\cdot\z}{x+\k\cdot\z}
{x+\k\cdot\z\choose \k}{y-\k\cdot\z\choose \n-\k} \nonumber\\
&\hskip -1.4cm={x+y\choose \n}+\sum_{i=1}^m
\sum_{\k=\0}^{\n}z_i {x-1+\k\cdot\z\choose \k-\e_i}{y-\k\cdot\z\choose \n-\k}.
\end{align*}
Then \eqref{eq:mh} follows from \eqref{eq:multi-2} by induction on
$|\n|$. However, I am unable to give a combinatorial proof of
Mohanty-Handa's identity.

Finally, we remark that a further generalization of
\eqref{eq:mh} was given by Chu \cite{Chu89} by using the following
generating functions due to Mohanty \cite{Mohanty66}:
\begin{align*}
&\sum_{\k\geq \0}\frac{x}{x+\k\cdot\z}
{x+\k\cdot\z\choose \k}u_1^{k_1}\cdots u_m^{k_m}
=v^{x}, \\
&\sum_{\k\geq \0}
{x+\k\cdot\z\choose \k}u_1^{k_1}\cdots u_m^{k_m}
=\frac{v^{x}}{1-\sum_{i=1}^m u_i z_i v^{z_i-1}},
\end{align*}
where $v$ satisfies the functional equation
$
\sum_{i=1}^m u_iv^{z_i}=v-1.
$


\begin{center}
{\small Department of Mathematics\\
East China Normal University\\
Shanghai 200062, People's Republic of China\\
E-mail: jwguo@math.ecnu.edu.cn\\
URL: http://math.ecnu.edu.cn/\textasciitilde{jwguo}}
\end{center}

\end{document}